\documentclass[12pt]{article}

\input{amssym.def}
\input{amssym.tex}

\newcommand{\bC}{{\bf C}}
\newcommand{\bT}{{\bf T}}
\newcommand{\bZ}{{\bf Z}}
\newcommand{\upi}{u_+^{-1}}
\newcommand{\umi}{u_-^{-1}}
\newcommand{\vpi}{v_+^{-1}}
\newcommand{\vmi}{v_-^{-1}}
\newcommand{\ai}{a^{-1}}
\newcommand{\vsk}{\vspace{1mm}}
\newcommand{\vsg}{\vspace{2mm}}

\newcommand{\Qnn}{Q_{-n}}

\begin{document}

\vspace*{2.5cm}
{\Large\bf
\begin{center}
On the determinant formulas by Borodin,\\
Okounkov, Baik, Deift, and Rains
\end{center}}

\vspace{2mm}
\begin{center}
{\large \bf A. B\"ottcher}
\end{center}

\vspace{2mm}
\begin{quote}

\renewcommand{\baselinestretch}{1.0}
\footnotesize
{We give alternative proofs to (block case versions of)
some formulas for Toeplitz and Fredholm determinants established
recently by the authors of the title.
Our proof of the Borodin-Okounkov formula is very short and direct.
The proof of the Baik-Deift-Rains formulas is based on standard
manipulations with Wiener-Hopf factorizations.}
\end{quote}

\vspace{8mm}
\noindent
{\large \bf 1. The formulas}

\vspace{4mm}
\noindent
Let $\bT$ be the complex unit circle and let $L^\infty:=L^\infty_{N
\times N}$ stand for the algebra of all $N \times N$ matrix
functions with entries in $L^\infty(\bT)$. Given $a \in L^\infty$,
we denote by $\{a_k\}_{k \in \bZ}$ the sequence of the Fourier
coefficients,
\[a_k=\frac{1}{2\pi}\int_0^{2\pi}a(e^{i
\theta})e^{-ik\theta}d\theta =\frac{1}{2\pi i}\int_{\bT}
a(z)z^{-k}\frac{dz}{z}.\]
The matrix function $a$ generates several structured (block)
matrices:
\begin{eqnarray*}
& & T(a) = (a_{j-k})_{j,k=0}^\infty \quad\mbox{(infinite block
Toeplitz),}\\[0.5ex]
& & T_n(a) = (a_{j-k})_{j,k=0}^{n-1} \quad\mbox{(finite block
Toeplitz),}\\[0.5ex]
& & H(a) = (a_{j+k+1})_{j,k=0}^\infty \quad\mbox{ (block
Hankel),}\\[0.5ex]
& & H(\tilde{a}) = (a_{-j-k-1})_{j,k=0}^\infty \quad\mbox{ (block
Hankel),}\\[0.5ex]
& & L(a) = (a_{j-k})_{j,k= -\infty}^\infty \quad\mbox{ (block
Laurent),}\\[0.5ex]
& & L(\tilde{a}) = (a_{k-j})_{j,k= -\infty}^\infty \quad\mbox{ (block
Laurent).}
\end{eqnarray*}
The matrices $T(a), H(a), H(\tilde{a})$ induce bounded operators
on $\ell^2(\bZ_+,\bC^N)$, and the matrices $L(a), L(\tilde{a})$
define bounded operators on $\ell^2(\bZ,\bC^N)$.

\pagebreak
\noindent
Let $\|\cdot\|$ be any matrix norm on $\bC^{N \times N}$.
We need the following classes of matrix functions:
\begin{eqnarray*}
& & W=\{a \in L^\infty: \sum_{n \in \bZ} \|a_n\| < \infty\} \quad
\mbox{(Wiener algebra),}\\[0.5ex]
& & K_1^1=\{a \in L^\infty: \sum_{n \in \bZ} (|n|+1)\|a_n\|<
\infty\} \quad\mbox{(weighted Wiener algebra),}\\[0.5ex]
& & K_2^{1/2}=\{a \in L^\infty: \sum_{n \in \bZ} (|n|+1)\|a_n\|^2<
\infty\} \quad\mbox{(Krein algebra),}\\[0.5ex]
& & H^\infty_\pm =\{a \in L^\infty: a_{\mp n} = 0 \;\:\mbox{for}\;\: n >0\}
\quad\mbox{(Hardy space).}
\end{eqnarray*}
Clearly, $K_1^1 \subset K_2^{1/2}$.
Given a subset $E$ of $L^\infty$, we say that a matrix function
$a \in L^\infty$ has a right (resp. left) canonical Wiener-Hopf
factorization in $E$ and write $a \in \Phi_r(E)$ (resp. $a \in
\Phi_l(E)$) if $a$ can be represented in the form $a=u_-u_+$
(resp. $a=v_+v_-$) with
\[u_-,v_-,\umi, v_-^{-1} \in E \cap H^\infty_-, \quad
u_+, v_+, \upi, v_+^{-1} \in E \cap H^\infty_+.\]
It is well known (see, e.g., \cite{BoSiBu},
\cite{WiII}) that if $a \in \Phi_r(L^\infty)$ then $T(a)$ is
invertible and $T^{-1}(a)=T(\upi)T(\umi)$ and that for $a \in
K_1^1$ (resp. $a \in W \cap K_2^{1/2}$) we have
\[a \in \Phi_r(K_1^1)\quad (\mbox{resp.}\;\:a \in
\Phi_r(K_2^{1/2})) \quad \Longleftrightarrow \quad
T(a) \;\: \mbox{is invertible}.\]
If $a \in K_1^1$ then $H(a)$ and $H(\tilde{a})$ are trace class
operators, and if $a \in K_2^{1/2}$, then  $H(a)$ and $H(\tilde{a})$
are Hilbert-Schmidt.

\vsk
We define the projections $P, Q, Q_n$ ($n \in \bZ$) on the space
$\ell^2(\bZ,\bC^N)$ by
\[(Px)_k=\left\{\begin{array}{cll} x_k & \mbox{for} & k \ge 0,\\
0 & \mbox{for} & k < 0,\end{array}\right. \qquad
(Qx)_k=\left\{\begin{array}{cll} 0 & \mbox{for} & k \ge 0,\\
x_k & \mbox{for} & k < 0,\end{array}\right.\]
\[(Q_nx)_k=\left\{\begin{array}{cll} x_k & \mbox{for} & k \ge n,\\
0 & \mbox{for} & k < n.\end{array}\right.\]
For $n \ge 1$, we let $P_n$ denote the projection on
$\ell^2(\bZ_+, \bC^N)$ given by
\[(P_nx)_k=\left\{\begin{array}{cll} x_k & \mbox{for} & 0 \le k \le n-1,\\
0 & \mbox{for} & k \ge n.\end{array}\right.\]
If $n \ge 0$, we can also think of $Q_n$ as an operator on
$\ell^2(\bZ_+,\bC^N)$. Note that the notation used here differs
from the one of \cite{BDR}, but that our notation is standard
in the Toeplitz business.

\vsk
On defining the flip operator $J$ on $\ell^2(\bZ,\bC^N)$ by
$(Jx)_k=x_{-k-1}$, we can write
\begin{equation}
T(a)=PL(a)P|{\rm Im}\,P, \quad
H(a)=PL(a)QJ|{\rm Im}\,P, \quad
H(\tilde{a}) =JQL(a)P|{\rm Im}\,P \label{1}
\end{equation}
Moreover, we may identify the operator $L(a)$ on
$\ell^2(\bZ,\bC^N)$ with the operator of multiplication
by $a$ on $L^2(\bT,\bC^N)$. Since $P, Q, J$ are also
naturally defined on the space $L^2(\bT,\bC^N)$, formulas (\ref{1})
enable us to interpret Toeplitz and Hankel operators
as operators on the Hardy space $H^2(\bT,\bC^N)$.

\vsk
For $a \in \Phi_l(L^\infty)$, the geometric mean $G(a)$
is defined by
$G(a)=(\det v_+)_0(\det v_-)_0,$
where $(\cdot)_k$ stands for the $k$th Fourier coefficient.
Thus, with an appropriately chosen logarithm,
\[G(a) = \exp (\log \det a)_0.\]
Let now $a \in \Phi_r(K_2^{1/2}) \cap \Phi_l(K_2^{1/2})$ and
let $a=u_-u_+$ and $a=v_+v_-$ be canonical Wiener-Hopf
factorizations. Put $b=v_-\upi$ and $c=\umi v_+$. Obviously,
$bc=I$. Using (\ref{1})
it is easily seen that
\begin{equation}
T(b)T(c)+H(b)H(\tilde{c})=I. \label{2}
\end{equation}
Since Hankel operators generated by matrix functions in
$K_2^{1/2}$ are Hilbert-Schmidt, the operator $H(b)H(\tilde{c})$
is in the trace class. From (\ref{2}) we infer that
$I-H(b)H(\tilde{c})$ is invertible. We put
\[E(a)=1/\det(I-H(b)H(\tilde{c})).\]
One can show (again see \cite{BoSiBu}, \cite{WiII}) that
$E(a)=\det T(a)T(\ai)$ and that in the scalar case ($N=1$)
we also have
\[E(a)=\exp \sum_{k=1}^\infty k(\log a)_k(\log a)_{-k}.\]

\vsg
\noindent
{\bf Theorem 1.1 (Borodin-Okounkov \`{a} la Widom).} {\it If
$a \in \Phi_r(K_2^{1/2}) \cap \Phi_l(K_2^{1/2})$ then
\begin{equation}
\det T_n(a) = G(a)^n E(a) \det(I-Q_nH(b)H(\tilde{c})Q_n) \label{3}
\end{equation}
for all $n \ge 1$.}

\vsg
In the scalar case, this beautiful theorem was established by
Borodin and Okounkov in \cite{BO}. It answered a question raised
by Its and Deift. The proof of \cite{BO} is rather complicated.
Three simpler proofs were subsequently found by Basor and Widom
\cite{BaWi} (who also extended the theorem to the block case)
and by the author \cite{Bot}. We here give still another proof,
which is very short and direct.

\vsk
Now suppose that $a \in \Phi_r(K_1^1) \cap \Phi_l(K_1^1)$
($ \subset \Phi_r(K_2^{1/2}) \cap \Phi_l(K_2^{1/2})$). Define $b$
and $c$ as above. We have
\[P-L(c)Q_nL(b)=(PL(c)-L(c)Q_n)L(b)
=(PL(c)Q-QL(c)P+L(c)(P-Q_n))L(b)\]
and since $PL(c)Q$ and $QL(c)P$ are trace class operators
(notice that $b,c \in K^1_1$) and the operator
$P-Q_n$ has finite rank, we see that $P-L(c)Q_nL(b)$ is
trace class.

\vsg
\noindent
{\bf Theorem 1.2 (Baik-Deift-Rains).} {\it If
$a \in \Phi_r(K_1^1) \cap \Phi_l(K_1^1)$ then
\begin{equation}
\det T_n(a) = G(a)^nE(a) 2^{-nN}\det(I+P-L(c)Q_nL(b)) \label{4}
\end{equation}
for all $n \ge 1$.}

\vsg
Clearly, to prove Theorem 1.2 it suffices to prove Theorem 1.1
and to verify that
\begin{equation}
\det (I+P-L(c)Q_nL(b))=2^{nN}\det (I-Q_nH(b)H(\tilde{c})Q_n)
\label{5}
\end{equation}
for all $n \ge 1$. By virtue of (\ref{1}),
\[\det (I-Q_nH(b)H(\tilde{c})Q_n) = \det (I-Q_nL(b)QL(c)Q_n)\]
for all $n \ge1$.
The right-hand side of the last equality makes sense for all
$n \in \bZ$. In fact, we have the following generalization of
(\ref{5}).

\vsg
\noindent
{\bf Theorem 1.3 (Baik-Deift-Rains).} {\it If $a \in
\Phi_r(K_1^1) \cap \Phi_l(K_1^1)$ then for all $n \in \bZ$,}
\begin{eqnarray}
& & \det (I+sP-sL(a)Q_nL(\ai))\nonumber\\
& & \quad =(1+s)^{nN}\det(I-s^2Q_nL(\ai)QL(a)Q_n)
\quad(s \neq -1) \label{6}\\
& & \quad = (1-s)^{-nN}\det(I-s^2(I-Q_n)L(\ai)PL(a)
(I-Q_n)) \quad (s \neq 1)\label{7}
\end{eqnarray}

\vsg
Theorems 1.2 and 1.3 are in \cite{BDR}. The proof given there
is as follows: the formulas are easily seen if some operator
that is no trace class operator were a trace class operator
and to save that insight the authors employ an approximation argument.
We here present a proof
that is a little more direct and uses Wiener-Hopf factorization.

\vsk
Theorem 1.1 is proved in Section 2, the proofs of
Theorems 1.2 and 1.3 are given in Section 3. In Section 4 we relax
the hypothesis of Theorem 1.3 to the requirement that $a$ be in $K_1^1$
and that $\det a$ have
no zeros on the unit circle, and in Section 5 we prove
a ``multi-interval'' version of Theorem 1.3.

\vspace{8mm}
\noindent
{\large \bf 2. Proof of the Borodin-Okounkov formula}

\vspace{4mm}
\noindent
If $K$ is an arbitrary trace class operator on $\ell^2(\bZ_+,\bC^N)$
and $I-K$ is invertible, then
\begin{equation}
\det P_n(I-K)^{-1}P_n = \frac{\det(I-Q_nKQ_n)}{\det(I-K)}.
\label{8}
\end{equation}
With $K$ replaced by $P_mKP_m$, this is Jacobi's theorem on the
principle $n \times n$ minor of the inverse of a (finite) matrix.
In the general case the identity follows from the fact
that $P_mKP_m$ converges to $K$ in the trace norm as $m \to \infty$.
For $K=H(b)H(\tilde{c})$ we obtain from (\ref{2}) that
\begin{eqnarray*}
& & P_n(I-K)^{-1}P_n =P_nT^{-1}(c)T^{-1}(b)P_n\\
& & =P_nT(v_+^{-1})T(u_-)T(u_+)T(v_-^{-1})P_n
=T_n(v_+^{-1})T_n(a)T_n(v_-^{-1}),
\end{eqnarray*}
and since $\det T_n(v_+^{-1})T_n(a)T_n(v_-^{-1}) = G(a)^{-n} \det T_n(a)$,
we get (\ref{3}) from (\ref{8}). \rule{2mm}{2mm}

\vspace{8mm}
\noindent
{\large \bf 3. Proof of the Baik-Deift-Rains formulas}

\vspace{4mm}
\noindent
In what follows we abbreviate $L(a)$ to $a$. Equivalently, we
may regard all operators on $L^2$ instead of $\ell^2$ and may
therefore think of $a$ as multiplication by $a$. Notice that
if $a \in K_1^1$ is invertible in $L^\infty$, then $\ai$ also
belongs to $K_1^1$.

\vsg
\noindent
{\bf Lemma 3.1.} {\it If $a$ and $\ai$ are in $K_1^1$ then
\[P-aQ_n\ai, \quad Q_n\ai QaQ_n, \quad (I-Q_n)\ai Pa(I-Q_n)\]
are trace class operators for all $n \in \bZ$.}

\pagebreak
\noindent
{\it Proof.} We have
\begin{eqnarray*}
& & P-aQ_n\ai =(Pa-aQ_n)\ai=(PaQ-QaP+a(P-Q_n))\ai,\\
& & Q_n\ai QaQ_n = -Q_n\ai QaP+Q_n\ai Qa(P-Q_n),\\
& & (I-Q_n)\ai Pa(I-Q_n)=(I-Q_n)\ai PaQ+(I-Q_n)\ai Pa(P-Q_n),
\end{eqnarray*}
and since $PaQ$ and $QaP$ are trace class and $P-Q_n$ has finite
rank, we arrive at the assertion. \rule{2mm}{2mm}

\vsk
We put \[f_n(s)=\det(I+sP-saQ_n\ai).\]

\vsg
\noindent
{\bf Proposition 3.2.} {\it If $a \in \Phi_r(K_1^1)$ and $n \ge 0$,
then}
\begin{equation}
f_n(s)=(1+s)^{nN}\det(I-s^2Q_n\ai QaQ_n). \label{21}
\end{equation}

\noindent
{\it Proof.} Let $a=u_-u_+$ be a right canonical Wiener-Hopf
factorization in $K^1_1$. Then
\[f_n(s)=\det(I+sP-su_-u_+Q_n\upi\umi)=
\det(I+s\umi P u_--su_+Q_n\upi),\]
and since $\umi Pu_-=P+Q\umi Pu_-P$ and $u_+Q_n=Pu_+Q_n$, we get
\[f_n(s)=\det(I+sQ\umi Pu_-P +sP-sPu_+Q_n\upi).\]
The operator $I+sQ\umi Pu_-P$ has the inverse $I-sQ\umi Pu_-P$
and its determinant is $1$. Hence,
\begin{eqnarray*}
f_n(s) & = & \det(I+(sP-sPu_+Q_n\upi)(I-sQ\umi Pu_-P))\\
& = & \det (I+sP-sPu_+Q_n\upi+s^2Pu_+Q_n\upi Q\umi Pu_-P).
\end{eqnarray*}
Because $\det(I+PA)=\det(I+PAP)$ and
\[Pu_-^{\pm 1}=Pu_-^{\pm 1}P,\quad
u_+^{\pm 1}P=Pu_+^{\pm 1}P,\quad
u_-^{\pm 1}Q=Qu_-^{\pm 1}Q,\quad
Qu_+^{\pm 1}=Qu_+^{\pm 1}Q,\]
it follows that
\begin{eqnarray*}
f_n(s) & = & \det (I+sP-sPu_+Q_n\upi P+s^2Pu_+Q_n\upi Q\umi
Pu_-P)\\
& = & \det (I+sP-sPu_+Q_n\upi P-s^2Pu_+Q_n\upi Q\umi Qu_-P)\\
& = & \det (I+sP-sQ_n-s^2P\upi Pu_+Q_n\upi Q\umi Qu_-Pu_+P)\\
& = & \det (I+sP-sQ_n-s^2Q_n\upi Q\umi Qu_-u_+P)\\
& = &  \det (I+sP-sQ_n-s^2Q_n\upi \umi Qu_-u_+P)\\
& = & \det (I+sP-sQ_n-s^2Q_n\ai QaP)\\
& = & \det (I+sP-sQ_n) \det(I-s^2Q_n\ai QaP)\\
& = & (1+s)^{nN} \det(I-s^2Q_n\ai QaQ_n) \quad \rule{2mm}{2mm}
\end{eqnarray*}

At this point we have proved formula (\ref{6}) for $n \ge 0$
and thus formula (\ref{5}) and Theorem 1.2. We are left
with switching from (\ref{6})
to (\ref{7}) and passing to negative $n$'s.

\pagebreak
\noindent
{\bf Proposition 3.3.} {\it If $a \in \Phi_l(K_1^1)$ and $n \ge 0$,
then}
\begin{equation}
f_{-n}(s)=(1-s)^{nN}\det(I-s^2(I-Q_{-n})\ai Pa(I-Q_{-n})).
\label{20}
\end{equation}

\noindent
{\it Proof.} We repeat the argument of the preceding proof, but
now we work with the left canonical Wiener-Hopf factorization
$a=v_+v_-$. We have
\begin{eqnarray*}
f_{-n}(s)& = & \det(I+sP-sa\Qnn \ai)\\
& = & \det(I-sQ+sa(I-\Qnn)\ai)\\
& = & \det(I-sQ+sv_+v_-(I-\Qnn)\vmi\vpi)\\
& = & \det(I-s\vpi Qv_++sv_-(I-\Qnn)\vmi)\\
& = & \det(I-sP\vpi Qv_+Q - sQ +sv_-(I-\Qnn)\vmi)\\
& = & \det(I+(-sQ +sQv_-(I-\Qnn)\vmi)(I+sP\vpi Qv_+Q))\\
& = & \det(I-sQ +sQv_-(I-\Qnn)\vmi Q
+s^2Qv_-(I-\Qnn)\vmi P\vpi Pv_+Q)\\
& = & \det(I-sQ +s(I-\Qnn)
-s^2Q\vmi Qv_-(I-\Qnn)\vmi P\vpi Pv_+Qv_-Q)\\
& = & \det(I-sQ +s(I-\Qnn)-s^2(I-\Qnn)\vmi \vpi Pv_+v_-Q)\\
& = & \det(I-sQ +s(I-\Qnn)-s^2(I-\Qnn)\ai PaQ)\\
& = & \det(I-sQ +s(I-\Qnn))\det(I-s^2(I-\Qnn)\ai PaQ)\\
& = & (1-s)^{nN} \det(I-s^2(I-Q_{-n})\ai Pa(I-Q_{-n})).
\quad \rule{2mm}{2mm}
\end{eqnarray*}

\vsg
\noindent
{\bf Lemma 3.4.} {\it If $a$ and $\ai$ are in $K^1_1$ and $n \in
\bZ$, then}
\begin{equation}
f_n(-s)f_n(s)=\det(I-s^2(I-Q_n)\ai Pa(I-Q_n))
\det (I-s^2Q_n\ai QaQ_n).
\label{9}
\end{equation}

\noindent
{\it Proof.} We have
\begin{eqnarray*}
& & (I-sP+saQ_n\ai)(I+sP-saQ_n\ai)\\
& & \quad = I-s^2P+s^2PaQ_n\ai P-s^2QaQ_n\ai Q\\
& & \quad = I-s^2Pa(I-Q_n)\ai P-s^2QaQ_n\ai Q.
\end{eqnarray*}
Taking determinants we obtain that
\begin{eqnarray*}
f_n(-s)f_n(s) & = & \det(I-s^2Pa(I-Q_n)\ai P) \det(I-s^2QaQ_n\ai Q)\\
& = & \det(I-s^2(I-Q_n)\ai Pa(I-Q_n))
\det (I-s^2Q_n\ai QaQ_n).
\quad \rule{2mm}{2mm}
\end{eqnarray*}

\vsg
\noindent
{\bf Proposition 3.5.} {\it If $a \in \Phi_l(K_1^1)$ and $n \ge 0$,
then
\begin{equation}
(1+s)^{nN}f_{-n}(s)= \det(I-s^2\Qnn\ai Qa\Qnn),
\label{18}
\end{equation}
and if  $a \in \Phi_r(K_1^1)$ and $n \ge 0$, then}
\begin{equation}
f_n(s)=(1-s)^{nN} \det(I-s^2(I-Q_n)\ai Pa(I-Q_n)).
\label{19}
\end{equation}

\noindent
{\it Proof.} Proposition 3.3 and Lemma 3.4 give
\begin{eqnarray*}
& & f_{-n}(s)(1+s)^{nN} \det(I-s^2(I-\Qnn)\ai Pa(I-\Qnn))
= f_{-n}(s)f_{-n}(-s)\\
& & \quad = \det(I-s^2(I-\Qnn)\ai Pa(I-\Qnn))
\det (I-s^2\Qnn\ai Qa\Qnn).
\end{eqnarray*}
Since $ \det(I-s^2(I-\Qnn)\ai Pa(I-\Qnn))\neq 0$ for sufficiently
small $s$, we get (\ref{18}) for these $s$ and then by analytic
continuation for all $s$. Analogously, using
Proposition 3.2 and Lemma 3.3 we get
\begin{eqnarray*}
& & f_n(s)(1-s)^{nN}\det(I-s^2Q_n\ai Q a Q_n) = f_n(s)f_n(-s)\\
& & \quad = \det(I-s^2(I-Q_n)\ai Pa(I-Q_n))
\det(I-s^2Q_n\ai QaQ_n),
\end{eqnarray*}
which implies (\ref{19}). \rule{2mm}{2mm}

\vsk
Theorem 1.3 is the union of Propositions 3.2, 3.4, and 3.5.

\vspace{8mm}
\noindent
{\large \bf 4. Non-invertible operators}

\vspace{4mm}
\noindent
The hypothesis of Theorem 1.3 is that $a$ be in $\Phi_r(K_1^1)\cap \Phi_l(K_1^1)$,
which is equivalent to the invertibility of both $T(a)$ and $T(\ai)$.
The theorem of this section, which is also from \cite{BDR},
relaxes this hypothesis essentially:
we only require that $T(a)$ be Fredholm
(which automatically implies that $T(\ai)$ is also
Fredholm). Notice that if
$a$ is continuous (and matrix functions in $K_1^1$ are continuous)
then $T(a)$ is a Fredholm operator if and only if $\det a$ has no zeros
on $\bT$. In that case the index of $T(a)$ is minus the winding number
of $\det a$, ${\rm Ind}\,T(a) = -{\rm wind}\,\det a$.

\vsg
\noindent
{\bf Lemma 4.1.} {\it If $a \in K_1^1$ and $T(a)$ is Fredholm
of index zero, then} (\ref{6}) {\it and} (\ref{7}) {\it are
valid.}

\vsg
\noindent
{\it Proof.} A theorem by Widom \cite{WiPe} tells us that there
exist a trigonometric polynomial $\varphi$ and a number
$\varrho>0$ such that $T(a +\varepsilon\varphi)$ is
invertible for all complex numbers $\varepsilon$ satisfying
$0<|\varepsilon|<\varrho$. Since $T(a+\varepsilon\varphi)$
is invertible, we conclude that
$a+\varepsilon\varphi \in \Phi_r(K_1^1)$.
Thus, (\ref{21}) and (\ref{19}) are true with $a$ replaced by
$a+\varepsilon\varphi$. From the proof of Lemma 3.1 we see
that
\begin{eqnarray*}
& & L(a+\varepsilon\varphi)Q_nL((a+\varepsilon\varphi)^{-1})
\to L(a)Q_nL(\ai),\\
& & Q_nL((a+\varepsilon\varphi)^{-1})QL(a+\varepsilon\varphi)Q_n
\to Q_nL(\ai)QL(a)Q_n
\end{eqnarray*}
in the trace norm as $\varepsilon \to 0$.
This gives (\ref{21}) and (\ref{19}). The proof of
formulas (\ref{20}) and (\ref{18}) is analogous. \rule{2mm}{2mm}

\vsg
\noindent
{\bf Lemma 4.2.} {\it If the scalar-valued
function $a \in K_1^1$ has no zeros on the unit circle and winding number
$w$ about the origin, then for all $n \in \bZ$,}
\begin{eqnarray}
& & \det (I+sP-sL(a)Q_nL(\ai))\nonumber\\
& & \quad =(1+s)^{n+w}\det(I-s^2Q_nL(\ai)QL(a)Q_n)
\quad(s \neq -1), \label{15}\\
& & \det (I+sP-sL(a)Q_nL(\ai))\nonumber\\
& & \quad =(1-s)^{-n-w}\det(I-s^2(I-Q_n)L(\ai)PL(a)
(I-Q_n)) \quad (s \neq 1).\label{16}
\end{eqnarray}

\noindent
{\it Proof.} Recall that $\chi_w$ is defined by $\chi_w(t)=t^w$.
We can write $a=\chi_wb$ with ${\rm wind}\,b=0$. The key
observation is that $\chi_wQ_n\chi_{-w}=Q_{n+w}$. Consequently,
\begin{eqnarray*}
& & \det(I+sP-saQ_n\ai)\\
& & \quad = \det(I+sP-sb\chi_wQ_n\chi_{-w}b^{-1})\\
& & \quad = \det (I+sP-sbQ_{n+w}b^{-1})\\
& & \quad = (1+s)^{n+w}\det(I-s^2Q_{n+w}b^{-1}QbQ_{n+w})
\quad \mbox{(by Theorem 1.3)}\\
& & \quad = (1+s)^{n+w}\det(I-s^2\chi_wQ_n\chi_{-w}b^{-1}Qb
\chi_wQ_n\chi_{-w})\\
& & \quad = (1+s)^{n+w}\det(I-s^2Q_n\ai QaQ_n),
\end{eqnarray*}
which is (\ref{15}). Analogously one can derive (\ref{16})
from Theorem 1.3. \rule{2mm}{2mm}

\vsg
\noindent
{\bf Theorem 4.3 (Baik-Deift-Rains).} {\it Let $a$ be an $N \times N$
matrix function in $K_1^1$ and suppose $\det a$ has no zeros on $\bT$.
Put $w = {\rm wind}\,\det a$. Then for all $n \in \bZ$,}
\begin{eqnarray}
& & \det (I+sP-sL(a)Q_nL(\ai))\nonumber\\
& & \quad =(1+s)^{nN+w}\det(I-s^2Q_nL(\ai)QL(a)Q_n)
\quad(s \neq -1), \label{25}\\
& & \det (I+sP-sL(a)Q_nL(\ai))\nonumber\\
& & \quad =(1-s)^{-nN-w}\det(I-s^2(I-Q_n)L(\ai)PL(a)
(I-Q_n)) \quad (s \neq 1).\label{26}
\end{eqnarray}

\noindent
{\it Proof} ({\it after Percy Deift}). We extend $a$ to an $(N+1)\times (N+1)$
matrix function $c$
by adding the $N+1,N+1$ entry $\chi_{-w}$:
\[c = \left(\begin{array}{cc} a & 0 \\ 0 & \chi_{-w} \end{array}
\right).\]
Since $T(c)$ is Fredholm of index zero, we deduce from Lemma 4.1 that
\begin{equation}
\det (I+sP-scQ_n c^{-1}) =(1+s)^{n(N+1)}\det(I-s^2Q_n c^{-1}Q c Q_n).
\label{27}
\end{equation}
Obviously,
\begin{eqnarray}
& & \det (I+sP-scQ_n c^{-1}) = \det (I+sP-saQ_n a^{-1})
\det (I+sP-s\chi_{-w}Q_n \chi_w),\label{28}\\
& & \det(I-s^2Q_n c^{-1}Q c Q_n) = \det(I-s^2Q_n a^{-1}Q a Q_n)
\det(I-s^2Q_n \chi_w Q \chi_{-w} Q_n). \label{30}
\end{eqnarray}
Lemma 4.2 implies that
\begin{equation}
\det (I+sP-s\chi_{-w}Q_n \chi_w) = (1+s)^{n-w}
\det(I-s^2Q_n \chi_w Q \chi_{-w} Q_n) \label{31}
\end{equation}
(which, by the way, can also be verified straightforwardly
in the particular case at hand).
Combining (\ref{27}), (\ref{28}), (\ref{30}), (\ref{31}) we arrive
at (\ref{25}). The proof of (\ref{26}) is analogous. \rule{2mm}{2mm}

\vspace{8mm}
\noindent
{\large \bf 5. The multi-interval case}

\vspace{4mm}
\noindent
The purpose of this section is to show that the argument
employed in Section 3 also works in the so-called multi-interval
case. The following theorem is again from \cite{BDR}.

\vsg
\noindent
{\bf Theorem 5 (Baik-Deift-Rains).} {\it Let $0 =n_0 \le n_1
\le \ldots \le n_k \le n_{k+1}=\infty$ be integers and let
$s_1, \ldots, s_k$ be complex numbers such that
$s_k-s_j \neq -1$ for all $j$. Put $s_0=0$. If $a \in
\Phi_r(K_1^1)$ then
\begin{eqnarray}
& & \hspace{-10mm}
\det\left(I+\sum_{j=1}^k(s_j-s_{j-1})(P-L(a)Q_{n_j}L(\ai))\right)
\nonumber\\
& & \hspace{-10mm}\quad
= \left(\,\prod_{j=0}^{k-1}(1\!+\!s_k\!-\!s_j)^{n_{j+1}-n_j}
\right) \det\left(I-\left(\,\sum_{j=1}^k\frac{s_ks_j}{1\!+\!s_k\!-\!s_j}
P_{[n_j,n_{j+1})}\right)L(\ai)QL(a)\right),
\label{17}
\end{eqnarray}
where $P_{[n_j,n_{j+1})}=Q_{n_j}-Q_{n_{j+1}}$ is the
projection onto the coordinates $l$ with $n_j \le l < n_{j+1}$.}

\pagebreak
\noindent
{\it Proof.} Proceeding exactly as in the proof of Proposition 3.2
we get
\begin{eqnarray*}
& & \det\left(I+\sum_{j=1}^k(s_j-s_{j-1})(P-aQ_{n_j}\ai)\right)\\
& & \quad = \det\left(I+\sum_{j=1}^k\left((s_j-s_{j-1})\umi P u_--
(s_j-s_{j-1})u_+Q_{n_j}\upi)\right)\right)\\
& & \quad = \det\left(I+\sum_{j=1}^k\left((s_j-s_{j-1})Q\umi P u_-P
+(s_j-s_{j-1})P-
(s_j-s_{j-1})Pu_+Q_{n_j}\upi)\right)\right)\\
& & \quad = \det\Bigg(I+\left(\sum_{j=1}^k(s_j-s_{j-1})P
- \sum_{j=1}^k(s_j-s_{j-1})Pu_+Q_{n_j}\upi \right)\\
& & \hspace{20mm} \times \left(I- \sum_{l=1}^k(s_l-s_{l-1})Q\umi
Pu_-P\right)\Bigg)\\
& & \quad = \det\Bigg(I+\sum_{j=1}^k(s_j-s_{j-1})P
-\sum_{j=1}^k(s_j-s_{j-1})Pu_+Q_{n_j}\upi\\
& & \hspace{20mm} +  \sum_{j,l}(s_j-s_{j-1})
(s_l-s_{l-1})Pu_+Q_{n_j}\upi Q\umi
Pu_-P\Bigg)\\
& & \quad = \det\Bigg(I+\sum_{j=1}^k(s_j-s_{j-1})P
-\sum_{j=1}^k(s_j-s_{j-1})Pu_+Q_{n_j}\upi P\\
& & \hspace{20mm} - \sum_{j,l}(s_j-s_{j-1})
(s_l-s_{l-1})Pu_+Q_{n_j}\upi Q\umi
Qu_-P\Bigg)\\
& & \quad = \det\Bigg(I+\sum_{j=1}^k(s_j-s_{j-1})P
-\sum_{j=1}^k(s_j-s_{j-1})Q_{n_j}\\
& & \hspace{20mm} -  \sum_{j,l}(s_j-s_{j-1})
(s_l-s_{l-1})Q_{n_j}\upi \umi Q u_-u_+P\Bigg).
\end{eqnarray*}
Clearly,
\[Q_{n_j}\upi \umi Q u_-u_+P=Q_{n_j}\ai QaP=:AP.\]
Since
\[\sum_{j=1}^k(s_j-s_{j-1})P
-\sum_{j=1}^k(s_j-s_{j-1})Q_{n_j}
= s_kP-\sum_{j=1}^ks_jP_{[n_j,n_{j+1})}
=\sum_{j=0}^{k-1}(s_k-s_j)P_{[n_j,n_{j+1})}\]
and
\[\sum_l(s_l-s_{l-1})=s_k,
\quad \sum_{j=1}^ks_k(s_j-s_{j-1})Q_{n_j}=
\sum_{j=1}^k s_k s_j P_{[n_j,n_{j+1})},\]
we obtain
\begin{eqnarray*}
& & \det\left(I+\sum_{j=1}^k(s_j-s_{j-1})(P-aQ_{n_j}\ai)\right)\\
& & \quad = \det\Bigg(I+\sum_{j=0}^{k-1}(s_j-s_k)P_{[n_j,n_{j+1})}
-\sum_{j=1}^ks_ks_jP_{[n_j,n_{j+1})}AP\Bigg)\\
& & \quad = \det\left(I+\sum_{j=0}^{k-1}(s_k-s_j)P_{[n_j,n_{j+1})}
\right)\\
& & \quad \hspace{10mm} \times \det\left(I-\left(\sum_{j=0}^k
\frac{1}{1+s_k-s_j}P_{[n_j,n_{j+1})}\right)\left(\,
\sum_{j=1}^k s_k s_j P_{[n_j,n_{j+1})}AP\right)\right)\\
& & \quad = \left(\,\prod_{j=0}^{k-1}(1\!+\!s_k\!-\!s_j)^{n_{j+1}-n_j}
\right) \det\left(I-\left(\,\sum_{j=1}^k\frac{s_ks_j}{1\!+\!s_k\!-\!s_j}
P_{[n_j,n_{j+1})}\right)A\right). \quad\rule{2mm}{2mm}
\end{eqnarray*}

\vsk
In \cite{BDR} it is also shown that if $a$ is a scalar-valued
function without zeros on the unit circle and with winding number
$w$, then (\ref{17})
is true with the additional factor
$(1+s_k)^w$
on the right-hand side. This can again be verified with the methods
developed here, but we stop at this point.

\vsk
\noindent
{\bf Acknowledgement.} I wish to thank Percy Deift and Harold Widom
for useful comments.

\vspace{4mm}
\noindent
\begin{minipage}[t]{7.5cm}
Fakult\" at f\" ur Mathematik \\
Technische Universit\"at Chemnitz \\09107 Chemnitz, Germany \\[0.5ex]
aboettch@mathematik.tu-chemnitz.de
\end{minipage}



\begin{thebibliography}{9}

\bibitem{BDR} J. Baik, P. Deift, and E.M. Rains: A Fredholm
determinant identity and the convergence of moments for random
Young tableaux. arXiv: math.CO/0012117.

\bibitem{BaWi} E.L. Basor and H. Widom: On a Toeplitz determinant
identity of Borodin and Oko\-unkov. {\it Integral Equations Operator
Theory} {\bf 37} (2000), 397-401.

\bibitem{BO} A. Borodin and A. Okounkov: A Fredholm determinant
formula for Toeplitz determinants. {\it Integral Equations Operator
Theory} {\bf 37} (2000), 386-396.

\bibitem{Bot} A. B\"ottcher:  One more proof of the Borodin-Okounkov
formula for Toeplitz determinants. arXiv: math.FA/0012200.

\bibitem{BoSiBu} A. B\"ottcher and B. Silbermann: {\it Analysis of
Toeplitz Operators.} Springer-Verlag, Berlin, Heidelberg, New
York 1990.

\bibitem{WiPe} H. Widom: Perturbing Fredholm operators to obtain
invertible operators. {\it J. Funct. Analysis} {\bf 20} (1975),
26-31.

\bibitem{WiII} H. Widom: Asymptotic behavior of block Toeplitz
matrices and determinants II. {\it Adv. Math.} {\bf 21} (1976),
1-29.

\end{thebibliography}
\end{document}